\documentclass[reqno,12pt,a4paper]{amsart}
\usepackage{amscd,amsfonts,amssymb,amsthm,latexsym}
\usepackage{srcltx}

\theoremstyle{plain}
\newtheorem{theorem}{Theorem}[section]
\newtheorem*{theorem*}{Theorem}
\newtheorem{corollary}{Corollary}[section]

\newtheorem{proposition}{Proposition}[section]

\theoremstyle{definition}
\newtheorem{definition}{Definition}[section]
\newtheorem*{definition*}{Definition}

\theoremstyle{remark}

\newtheorem*{remark*}{Remark}

\numberwithin{equation}{section}

\begin{document}
\raggedbottom 

\title[Sidon-type inequalities]{Sidon-type inequalities and the space of quasi-continuous functions}

\author{Artyom Radomskii}

\begin{abstract}
We discuss some different results on Sidon-type inequalities and on the space of quasi-continuous functions.
\end{abstract}

 \address{Steklov Mathematical Institute of Russian Academy of Sciences\\
 8 Gubkina St., Moscow 119991, Russia}

\thanks{This work was performed at the Steklov International Mathematical Center
and supported by the Ministry of Science and Higher Education of the
Russian Federation (agreement no. 075-15-2019-1614).}

\email{artyom.radomskii@mi-ras.ru}

\maketitle

\begin{flushright}
\emph{To Boris Sergeevich Kashin on his 70th birthday}
\end{flushright}

\section{Introduction}

 In this paper we discuss some different results on Sidon-type inequalities and on the space of quasi-continuous functions. Let us introduce some notation. The sets of positive integers, integers, and real numbers are denoted by $\mathbb{N}$, $\mathbb{Z}$, and $\mathbb{R}$ respectively. Given $f\in L^{p}(0,2\pi)$, we put
\begin{align*}
\|f\|_{p}&=\left(\int_{0}^{2\pi}|f(x)|^{p}\, dx\right)^{1/p}\ \ \text{for $1\leq p <\infty$},\\
\|f\|_{\infty}&= \text{ess} \sup_{[0,2\pi]} |f(x)|\ \ \ \ \,\ \ \ \ \ \ \ \text{for $p=\infty$}.
\end{align*} Here
\[
\text{ess} \sup_{[0,2\pi]} |f(x)|= \inf \Bigl\{ C>0:\ \textup{mes}\bigl\{x\in [0, 2\pi]:\ |f(x)| > C\bigr\} = 0\Bigr\},
\]where $\textup{mes}$ denotes the Lebesgue measure. We observe that if $f$ is a continuous function on $[0,2\pi]$, then
\[
\|f\|_{\infty}=\max_{[0,2\pi]}|f(x)|.
\]

If $x$ is a real number, then $[x]$ denotes its integral part, and $\lceil x\rceil$ is the smallest integer $n$ such that $n\geq x$. We put $\log_{a}x:=\ln x/\ln a$. For a nonzero trigonometric polynomial $T(x)$, its exact order will be denoted by $\deg (T)$. For a real number $r\geq 0$, let $\mathrm{T}(r)$ denote the space of all real trigonometric polynomials of the form
\[
t(x)=A+\sum_{k=1}^{[r]}(a_k\cos kx + b_k\sin kx)
\](by definition, the sum $\sum_{k=1}^{0}$ is put to be zero). For a positive integer $n$, we denote by $E_n$ the space of all real trigonometric polynomials of the form
\[
t(x)=\sum_{k=2^{n-1}}^{2^{n}-1} (a_k\cos kx+b_k\sin kx).
\]By definition, put $E_0=\mathrm{T}(0)$. We write $\dim V$ for the dimension of a finite dimensional vector space $V$ over $\mathbb{R}$. It is easy to see that $\mathrm{T}(n)$ and $E_n$ are finite dimensional vector spaces and that $\dim \mathrm{T}(n)=2n+1$ and $\dim E_n=2^n$. It is also evident that
\[
 \mathrm{T}(2^n - 1)= E_0\oplus E_1 \oplus\cdots\oplus E_n.
 \]

Let $\lambda$ be a real number with $\lambda>1$. We denote by $\Lambda(\lambda)$ the class of sequences $U=\{n_k\}_{k=1}^{\infty}$ of positive integers such that $n_{k+1}/n_k\geq \lambda$, $k=1, 2, \ldots .$ Let $\Lambda$ stand for the class of all lacunary sequences $U$, i.e.,
\[
\Lambda=\bigcup\limits_{\lambda>1} \Lambda(\lambda).
\]Finally, we introduce the class $\Lambda_\sigma$ of all increasing sequences $U$ of positive integers that admit splitting into finitely many lacunary sequences. We observe that if $U\in \Lambda_{\sigma}$, then $U$ is an increasing sequence of positive integers, and it can be split into finitely many sequences $U^{(j)}\in \Lambda(\mu)$, where $\mu>1$ is any number prescribed beforehand.

We recall that functions $f, g\in L^{2}(0, 2\pi)$ are \emph{orthogonal} if
\[
\int_{0}^{2\pi} f(x)\overline{g(x)}\, dx = 0
\](here the overbar denotes complex conjugation). A system of functions $\{\varphi_{n}(x)\}_{n=1}^{\infty}\subset L^{2}(0,2\pi)$ is \emph{orthonormal} if
\[
\int_{0}^{2\pi} \varphi_{n}(x)\overline{\varphi_{m}(x)}\, dx =
\begin{cases}
1 &\text{if $n=m$,}\\
0 &\text{if $n\neq m$},
\end{cases}\quad n, m = 1, 2, \ldots.
\]

Given $f\in L^{1}(0,2\pi),$ we consider the Fourier coefficients
\begin{gather*}
a_{k}(f)=\frac{1}{\pi}\int_{0}^{2\pi}f(x)\cos kx\,dx,\quad b_{k}(f)=\frac{1}{\pi}\int_{0}^{2\pi}f(x)\sin kx\,dx\\
\qquad\qquad\qquad\qquad\qquad\qquad\qquad\qquad\qquad\qquad\qquad(k=0, 1, 2, \ldots).
\end{gather*}

\section{Sidon-type inequalities for trigonometric polynomials}

In 1927 Sidon proved the following result.

\begin{theorem}[Sidon \cite{Sidon1}]\label{Th.Sidon}
 Let $\{n_k\}_{k=1}^{\infty}$ be a sequence of positive integers such that
\[
\frac{n_{k+1}}{n_k}\geq \lambda>1,\quad k=1, 2,\ldots.
\]If a trigonometric series
\[
\sum_{k=1}^{\infty}(\alpha_{k}\cos n_k x + \beta_k \sin n_k x)
\]is the Fourier series of a bounded measurable function $f(x)$, then
\[
\sum_{k=1}^{\infty}(|\alpha_k|+|\beta_k|)<\infty.
\]
\end{theorem} Sidon's method of proof was based on application of the Riesz products, which became an important tool in the theory of trigonometric and general orthogonal series. For the first time, these products arose in F. Riesz's paper \cite{RieszF1}. Also in \cite{Sidon1}, Sidon observed that Theorem \ref{Th.Sidon} remains valid in the case where $f(x)$ is only bounded from one side, i.e., if $f(x)\leq M$ or $f(x)\geq -M$. In fact, the proof of Theorem \ref{Th.Sidon} in \cite{Sidon1} implies the following estimate:
\[
\sum_{k=1}^{\infty}(|\alpha_k|+|\beta_k|)\leq C(\lambda)\|f\|_{\infty},
\]where $C(\lambda)>0$ is a constant depending only on $\lambda$. In particular, the next result is true.
\begin{theorem}
Let $\{n_k\}_{k=1}^{\infty}$ be a sequence of positive integers such that
\[
\frac{n_{k+1}}{n_k}\geq \lambda>1,\quad k=1, 2,\ldots.
\]Let $m$ be a positive integer. Let $\alpha_k, \beta_k\in \mathbb{R},$ $k=1,\ldots, m,$ and let
\[
f(x)=\sum_{k=1}^{m}(\alpha_k\cos n_k x+ \beta_k\sin n_k x).
\]Then
\begin{equation}\label{S1:SidonIneq}
\|f\|_{\infty}\geq c(\lambda)\sum_{k=1}^{m}(|\alpha_k|+|\beta_k|),
\end{equation}where $c(\lambda)>0$ is a constant depending only on $\lambda$.
\end{theorem}

A direction of refinement of Theorem \ref{Th.Sidon} was related to relaxing its assumptions concerning lacunarity. In \cite{Sidon2}, Sidon himself carried his theorem over to the case where $U=\{n_k\}_{k=1}^{\infty}$ can be split into finitely many lacunary sequences, i.e. $U\in \Lambda_{\sigma}$. The further developments in this issue were done in the papers \cite{Stechkin}, \cite{Pisier_1} -- \cite{Pisier_3}, and others.

In 1998, Kashin and Temlyakov \cite{Kashin_Temlyakov1,Kashin_Temlyakov2} started the study of another direction of refining the Sidon theorem. In connection with estimates for the entropy numbers of some classes of functions of small smoothness, they explored the question about the possible generalizations of the Sidon inequality \eqref{S1:SidonIneq} where $\alpha_k \cos n_k x$ is replaced by $p_k (x)\cos n_k x$, $p_k(x)$ being a trigonometric polynomial.

\begin{theorem}[Kashin and Temlyakov \cite{Kashin_Temlyakov1,Kashin_Temlyakov2}]\label{Th.KashTeml}
Let $l$ be a positive integer. Let $p_k\in \mathrm{T}(2^l),$ $k=l+1,\ldots, 2l$, and let
\[
f(x)=\sum_{k=l+1}^{2l}p_{k}(x)\cos 4^k x.
\]Then
\begin{equation}\label{KashTemIneq}
\|f\|_{\infty}\geq c\sum_{k=l+1}^{2l}\|p_k\|_1,
\end{equation}where $c>0$ is an absolute constant.
\end{theorem}

In \cite{Radomskii.Algebra.An} we carried the result of Theorem \ref{Th.KashTeml} to sequences of class $\Lambda_\sigma$ and relaxed the conditions imposed on the degrees of the trigonometric polynomials $p_k(x)$. This refinement was achieved via application of a new method of proof, based on a refinement of the Riesz products and an estimate for the number of solutions of certain Diophantine equations.
\begin{theorem}[Radomskii \cite{Radomskii.Algebra.An}]\label{S1:T.Rad.Alg}
Let $\varepsilon$ and $B$ be real numbers with $\varepsilon\in (0,1)$ and $B\geq 1$, and let $U=\{n_k\}_{k=1}^{\infty}$ be an increasing sequence of positive integers such that $U$ can be split into $d$ sequences $U^{(j)}\in \Lambda (\lceil 7\sqrt{B}\rceil),$ $j=1,\ldots,d.$ Let $l$ and $m$ be positive integers with $m\geq l+2$. Let $p_k, q_k\in \mathrm{T}(r_k),$ $k=l,\ldots, m,$ where
\begin{align*}
r_l &= \min\left(\frac{n_{l+1}-n_{l}}{2(1+\varepsilon)}, \frac{n_l}{1+\varepsilon}\right),\\
 r_k &= \min\left(\frac{n_{k}-n_{k-1}}{2(1+\varepsilon)}, \frac{n_{k+1}-n_{k}}{2(1+\varepsilon)}, B n_{l}\right),\ k=l+1,\ldots, m-1,\\
r_m &= \min\left(\frac{n_{m}-n_{m-1}}{2(1+\varepsilon)}, B n_{l}\right),
\end{align*} and let
\[
f(x)=\sum_{k=l}^{m} \big(p_{k}(x)\cos n_k x + q_k(x)\sin n_k x\big).
\]Then
\[
\|f\|_{\infty}\geq \frac{c}{d^{2}\cdot\ln^{2}(1+1/\varepsilon)}\,\sum_{k=l}^{m}(\|p_k\|_{1}+\|q_k\|_{1}),
\]where $c>0$ is an absolute constant.
\end{theorem}
\begin{corollary}\label{S1:Cor}
Let $\varepsilon$ be a real number with $\varepsilon\in (0,1)$, and let $\{n_{k}\}_{k=1}^{\infty}$ be a sequence of positive integers such that
\[
\frac{n_{k+1}}{n_k}\geq \lambda> 1,\quad k=1, 2, \ldots .
\]We put
\[
\gamma = \min \left(\frac{\lambda -1}{2(1+\varepsilon)}, \frac{1}{1+\varepsilon}\right).
\]Let $l$ and $m$ be positive integers with $m\geq l+2$. Let $p_k, q_k\in \mathrm{T}(\gamma n_l),$ $k=l,\ldots, m,$ and let
\[
f(x)=\sum_{k=l}^{m} \big(p_{k}(x)\cos n_k x + q_k(x)\sin n_k x\big).
\]Then
\[
\|f\|_{\infty}\geq \frac{c}{\lceil\ln 7/\ln \lambda\rceil^{2}\cdot\ln^{2}(1+1/\varepsilon)}\,\sum_{k=l}^{m}(\|p_k\|_{1}+\|q_k\|_{1}),
\]where $c>0$ is an absolute constant.
\end{corollary}
\textsc{Proof of Corollary \ref{S1:Cor}.} We take $B=1$, and put
\begin{gather*}
d=\lceil \ln 7/ \ln \lambda\rceil,\\
U^{(j)}= \{n_{j+(k-1)d}\}_{k=1}^{\infty},\ j=1,\ldots, d.
\end{gather*}For $1\leq j \leq d$ and $k\in \mathbb{N}$, we have
\[
\frac{n_{j+kd}}{n_{j+(k-1)d}}=\frac{n_{j+kd}}{n_{j+kd-1}}\cdots\frac{n_{j+kd-d+1}}{n_{j+kd-d}}\geq \lambda^{d}\geq 7.
\]We see that $U=\{n_k\}_{k=1}^{\infty}$ is an increasing sequence of positive integers and $U$ can be split into $d$ sequences $U^{(j)}\in \Lambda (7),$ $j=1, \ldots, d.$ It is easy to see that $r_k\geq \gamma n_l,$ $k=l,\ldots, m.$ Therefore $p_k, q_k\in \mathrm{T}(r_k),$ $k=l,\ldots, m.$ From Theorem \ref{S1:T.Rad.Alg} we have
\begin{align*}
\|f\|_{\infty}&\geq \frac{c}{d^{2}\cdot\ln^{2}(1+1/\varepsilon)}\,\sum_{k=l}^{m}(\|p_k\|_{1}+\|q_k\|_{1})=\\
&=\frac{c}{\lceil\ln 7/\ln \lambda\rceil^{2}\cdot\ln^{2}(1+1/\varepsilon)}\,\sum_{k=l}^{m}(\|p_k\|_{1}+\|q_k\|_{1}),
\end{align*}where $c>0$ is an absolute constant. Corollary \ref{S1:Cor} is proved.

\section{Some trigonometric polynomials with extremely small uniform norm}

We first note the following fact. Let $\alpha$ be a positive real number and $\{t_{n}(x)\}_{n=1}^{\infty}$ be a sequence of orthogonal trigonometric polynomials such that $\|t_n\|_1 \geq \alpha,$ $n=1, 2, \ldots .$ Then
\begin{align*}
\biggl\|\sum_{j=1}^{n}t_{j}\biggr\|_{\infty}&\geq\frac{1}{\sqrt{2\pi}}\biggl\|\sum_{j=1}^{n}t_{j}\biggr\|_{2}=
\frac{1}{\sqrt{2\pi}}\biggl(\sum_{j=1}^{n}\|t_{j}\|_{2}^{2}\biggr)^{1/2}\geq\\
&\geq\frac{1}{\sqrt{2\pi}}\biggl(\frac{1}{2\pi}\sum_{j=1}^{n}\|t_{j}\|_{1}^{2}\biggr)^{1/2}\geq\frac{\alpha}{2\pi}\sqrt{n},\quad n=1,2,\ldots.
\end{align*}Thus, under the assumptions above, the uniform norm of $\sum_{j=1}^{n}t_j$ cannot have order of growth smaller than $\sqrt{n}$.

In 1997 Grigor'ev proved the following result.

\begin{theorem}[Grigor'ev \cite{Grig}]\label{S3:Th.Grig}
   There is a sequence $\{t_n(x)\}_{n=1}^{\infty}$ of complex trigonometric polynomials such that
  \begin{gather*}
  t_n(x)=\sum\limits_{2^n\leq |j|<2^{n+1}} c_j e^{i j x},\quad \|t_{n}\|_{1}\geq \frac{\pi}{4},\quad \|t_n\|_{\infty}\leq 6,\quad n=1,2,\ldots,\\
  \biggl\|\sum_{k=1}^{n}t_k \biggr\|_{\infty}\leq A\sqrt{n},\quad n=1,2,\ldots,
  \end{gather*}
  where $A>0$ is an absolute constant.
  \end{theorem}The technique proposed in \cite{Grig} can be called the \emph{pseudo stopping time method}. In \cite{Radomskii.Pos.Str} we introduced a refinement of the Grigor'ev method for studying the possibility of strengthening Sidon-type inequalities. This refinement allowed us to construct trigonometric polynomials as in Theorem \ref{S3:Th.Grig} but with a stronger restriction on the spectrum. The result in \cite{Radomskii.Pos.Str} showed that in Corollary \ref{S1:Cor} (for $n_k=2^k$) it is impossible to replace the condition $p_k, q_k\in \mathrm{T}(\gamma 2^l)$ by $p_k, q_k\in \mathrm{T}(2^{k-k^{\theta}})$ for any $\theta\in (0,1).$ The further investigation of these issues was continued in \cite{GrigRad2015}, \cite{Radomskii.Izv}. The best result is as follows.
  \begin{theorem}[Radomskii \cite{Radomskii.Izv}]\label{S3:Th.Rad.Izv}
  Let $\lambda$ be a real number with $\lambda >1$, $f:\ [1, +\infty)\rightarrow \mathbb{R}$ be a positive non-increasing function, and $\{m_n\}_{n=1}^{\infty}$ be a sequence of positive integers such that
  \[
\frac{m_{n+1}}{m_n}\geq \lambda^{f(n)}=:\lambda_n,\quad n=1, 2, \ldots.
\]Also, let a sequence $\{d_n\}_{n=1}^{\infty}$ of positive integers satisfy the condition $1\leq d_n\leq m_{n+1}-m_n,$ $n=1,2,\ldots .$ Then there is a sequence $\{t_n(x)\}_{n=1}^{\infty}$ of complex trigonometric polynomials such that
\begin{gather*}
t_{n}(x)=\sum_{m_n\leq s < m_n+d_n} c_s e^{i s x},\quad \|t_n\|_1\geq \frac{5\pi}{8},\quad \|t_n\|_{\infty}\leq 6,\quad n=1, 2,\ldots,\\
\biggl\|\sum_{j=1}^{n}t_{j}\biggr\|_{\infty}\leq 277 + \alpha\sqrt{n}+180\max\limits_{1\leq j\leq n}\log_{\lambda_j}\max\biggl(\frac{m_j}{d_j}, \frac{1}{\ln \lambda_j}, 1\biggr),\quad n=1, 2,\ldots,
\end{gather*}where $\alpha>0$ is an absolute constant.
  \end{theorem}With minimal changes in the proof, one can prove the modified version of the theorem with $t_{j}(x) = p_j(x) \cos m_j x$, where the $p_j\in \mathrm{T}(d_j)$. From Theorem \ref{S3:Th.Rad.Izv} we obtain
  \begin{corollary}[Radomskii \cite{Radomskii.Izv}]
  Let $\lambda$, $\varepsilon$, and $\theta$ be real numbers with $\lambda>1,$ $0\leq \varepsilon<1,$ $0<\theta < 1-\varepsilon,$ and let
  \[
m_n=\big[\lambda^{n^{1-\varepsilon}}\big],\quad d_n=\big[\lambda^{n^{1-\varepsilon}-n^{\theta}}\big],\qquad n=1,2,\ldots.
\]Then there is a sequence $\{t_n(x)\}_{n=1}^{\infty}$ of complex trigonometric polynomials such that
\begin{gather*}
t_{n}(x)=\sum_{m_n\leq s < m_n+d_n} c_s e^{i s x},\quad \|t_n\|_1\geq \frac{5\pi}{8},\quad \|t_n\|_{\infty}\leq 6,\quad n=1, 2,\ldots,\\
\biggl\|\sum_{j=1}^{n}t_{j}\biggr\|_{\infty}\leq C(\lambda, \varepsilon, \theta)\,n^{\max(1/2,\,\theta+\varepsilon)},\quad n=1, 2,\ldots,
\end{gather*}where $C(\lambda,\varepsilon,\theta)>0$ is a constant depending only on $\lambda$, $\varepsilon$, and  $\theta$.
  \end{corollary}Also in \cite{Radomskii.Izv} we proved the following result.
  \begin{theorem}[Radomskii \cite{Radomskii.Izv}]\label{Th.Rad.Poly.Ln}
  Let $\varepsilon$ be a real number with $\varepsilon\in (0,1)$, $\{k_n\}_{n=1}^{\infty}$ be a sequence of positive integers with  $k_1< \cdots< k_n<\cdots$, and $L_n$ be a subspace of $E_{k_n}$ such that $\dim L_n \geq \varepsilon \dim E_{k_n},$ $n=1, 2, \ldots .$ Then there is a sequence $\{t_n(x)\}_{n=1}^{\infty}$ of trigonometric polynomials such that
  \begin{gather*}
t_n\in L_n,\quad \|t_n\|_{\infty}\leq 1,\quad \|t_n\|_1\geq c\cdot\varepsilon,\quad n=1,2,\ldots,\\
\biggl\|\sum_{j=1}^{n} t_j\biggr\|_{\infty}\leq \frac{a}{\sqrt{\varepsilon}}\sqrt{n}+\frac{b}{\varepsilon^{3/2}},\quad n=1,2,\ldots,
\end{gather*}where $a$, $b$, $c$ are positive absolute constants.
  \end{theorem}

  \section{The Rademacher functions}

  In this section we give some well-known facts on the Rademacher functions.
  \begin{definition}
  A set of real measurable functions $\{f_n(x)\}_{n=1}^{N}$ with domain $(0,1)$ is a \emph{set of independent functions} if for every interval $I_{n},$ $n=1,\ldots, N,$ of the real line the following condition is satisfied:
  \begin{gather*}
  \textup{mes}\{x\in (0,1): f_{n}(x)\in I_{n},\ n=1,\ldots, N\}=\\
  =\prod_{n=1}^{N}\textup{mes}\{x\in (0,1): f_{n}(x)\in I_{n}\}.
  \end{gather*}

  An infinite sequence of functions $\{f_{n}(x)\}_{n=1}^{\infty}$ is a \emph{system of independent functions} if the set $\{f_{n}(x)\}_{n=1}^{N}$ is a set of independent functions for every $N=1, 2, \ldots .$
  \end{definition}

  Systems of independent functions have many interesting properties. For example (see \mbox{\cite[Chapter 2]{KashSaak}}), for every set of independent functions $\{f_{n}(x)\}_{n=1}^{N}$ with $f_n\in L^{1}(0,1),$ $n=1,\ldots, N,$ the function $\prod_{n=1}^{N}f_{n}(x)$ also belongs to $L^{1}(0,1)$ and
  \[
 \int_{0}^{1}\prod\limits_{n=1}^{N} f_{n}(x)\, dx= \prod\limits_{n=1}^{N} \int_{0}^{1} f_{n}(x)\,dx.
 \]
  One of the most important systems of independent functions is the Rademacher system $\{r_{n}(x)\}_{n=1}^{\infty}$.
  \begin{definition}
  For $n=1, 2, \ldots ,$ the $n$th Rademacher function is defined by
  \begin{align*}
 r_{n}(x)&=
 \begin{cases}
 1   &\text{for $x\in((i-1)/2^{n}, i/2^{n}),\  \text{$i$ odd}$}\\
 -1  &\text{for $x\in((i-1)/2^{n}, i/2^{n}),\ \text{$i$ even}$}
 \end{cases}\\
 &\qquad\qquad\qquad\qquad\qquad\qquad\qquad\ (i=1,\ldots, 2^n).
 \end{align*}
  \end{definition}In addition, it will be convenient to suppose in what follows that $r_{0}(x)=1$ for $x\in (0,1)$ and that $r_{n}(i/2^{n})=0$ for $i=0, 1,\ldots, 2^n ;$ $n=0, 1, 2,\ldots .$

  The Rademacher functions were introduced in 1922 in the paper of Rademacher \cite{Rademacher}. The important properties of these functions are given in the following proposition (see, for example, \mbox{\cite[Chapter 2]{KashSaak}}).

  \begin{proposition}\label{S4:Rademacher.Prop}
  The following statements hold.\\
   1) The functions $\{r_{n}(x)\}_{n=0}^{\infty},$ $x\in [0,1],$ form a system of independent functions.\\
   2) The system of functions $\{r_{n}(x)\}_{n=0}^{\infty},$ $x\in [0,1],$ is orthonormal.\\
   3) Let $a_1,\ldots,a_n$ and $\lambda$ be real numbers, $\lambda\geq 0$. Then
 \begin{align*}
 &\text{i) } \textup{mes}\biggl\{x\in[0,1]:\ \biggl|\sum_{j=1}^{n}a_j r_{j}(x)\biggr| > \lambda \biggl(\sum_{j=1}^{n} a_{j}^{2}\biggr)^{1/2}\biggr\}\leq 2 e^{-\lambda^{2}/2};\\
 &\text{ii) } \int_{0}^{1} \biggl|\sum_{j=1}^{n}a_j r_{j}(x)\biggr|\, dx\geq c_1 \biggl(\sum_{j=1}^{n} a_{j}^{2}\biggr)^{1/2};\\
 &\text{iii) } \int_{0}^{1}\max\limits_{1\leq k \leq n} \biggl|\sum_{j=1}^{k}a_j r_{j}(x)\biggr|\,dx\leq
 c_2 \biggl(\sum_{j=1}^{n} a_{j}^{2}\biggr)^{1/2},
 \end{align*} where $c_1>0$ and $c_2>0$ are absolute constants.
  \end{proposition}

\section{The space of quasi-continuous functions}

For $f\in L^1(0, 2\pi)$ with Fourier series $f\sim \sum_{n=0}^{\infty}\delta_{n}(f,x)$, where
\begin{gather*}
 \delta_{0}(f)=a_0 (f)/2,\\
 \delta_{n}(f, x)=\sum_{k=2^{n-1}}^{2^n - 1} \bigl(a_k (f)\cos kx + b_k(f) \sin kx\bigr),\quad n=1,2,\ldots,
 \end{gather*}we define its $\mathrm{QC}$-norm by
  \begin{equation}\label{S5:QC.Norm.DEF}
 \|f\|_{\mathrm{QC}}\equiv \int_{0}^{1} \biggl\|\sum_{n=0}^{\infty}r_n (\omega)\delta_{n}(f,\cdot)\biggr\|_{\infty}\, d\omega.
 \end{equation}Here $r_{n}(\omega)$ is the $n$th Rademacher function. By the \emph{space of quasi-continuous functions} we mean the closure of the set of trigonometric polynomials with respect to the norm \eqref{S5:QC.Norm.DEF}. The space of quasi-continuous functions and the $\mathrm{QC}$-norm were introduced by Kashin and Temlyakov in \cite{Kashin_Temlyakov1} and \cite{Kashin_Temlyakov2}. They proved the following result.
 \begin{theorem}[Kashin and Temlyakov \cite{Kashin_Temlyakov1,Kashin_Temlyakov2}]\label{S5:Th.KashTem}
 Let $n$ be a positive integer. Let $t_j\in E_j,$ $j=0,\ldots, n,$ and let $f(x)=\sum_{j=0}^{n} t_{j}(x).$ Then
 \[
 \|f\|_{\mathrm{QC}}\geq \frac{1}{48\pi}\sum_{j=0}^{n}\|t_j\|_1.
 \]
 \end{theorem}In connection with problems of approximation theory (see \cite{Kashin_Temlyakov2} for details) the question of the connection between the $\mathrm{C}$- and $\mathrm{QC}$-norms is of interest. Theorem \ref{S3:Th.Grig} and Theorem \ref{S5:Th.KashTem} imply
 \begin{theorem}[\cite{Kashin_Temlyakov2}]Let $n$ be a positive integer. Then
 \[
 \sup\limits_{t\in \mathrm{T}(2^n)} \frac{\|t\|_{\mathrm{QC}}}{\|t\|_{\infty}}\geq c_1 \sqrt{n},
 \]where $c_1>0$ is an absolute constant.
 \end{theorem} Oskolkov (see \cite{Kashin_Temlyakov2}) proved
 \begin{theorem}\label{S5:Th.Oskolkov}
 Let $n$ be a positive integer. Then
 \[
\sup\limits_{t\in \mathrm{T}(2^n)} \frac{\|t\|_{\infty}}{\|t\|_{\mathrm{QC}}}\geq c_2 \sqrt{n},
\]where $c_2>0$ is an absolute constant.
\end{theorem}
\textsc{Sketch of proof of Theorem \ref{S5:Th.Oskolkov}.} Let us consider
\[
t(x)=\sum_{j=1}^{n}\frac{1}{2^{j-1}}\sum_{k=2^{j-1}}^{2^{j}-1}\cos kx=\sum_{j=1}^{n}\delta_{j}(x).
\]It is clear that $t\in\mathrm{T}(2^n)$ and $\|t\|_{\infty}= t(0)=n.$ We put
\[
t_{\omega}(x)=\sum_{j=1}^{n}r_{j}(\omega) \delta_{j}(x),\qquad
g_{\omega}(x)=\sum_{j=1}^{n} r_{j}(\omega)\chi_{[-\pi 2^{-j}, \pi 2^{-j}]}(x),
\]where $\chi_{E}(x)$ is a characteristic function of a set $E$. Using the inequality $|1-\cos x|\leq |x|,$ one can show that 
 \[
 \|t_{\omega}(x) - g_{\omega}(x)\|_{\infty}\leq c
 \]for any $\omega \in [0,1]$, where $c>0$ is an absolute constant. Applying Proposition \ref{S4:Rademacher.Prop} (3, iii), we have
\[
 \int_{0}^{1}\|g_{\omega}(x)\|_{\infty}\, d\omega= \int_{0}^{1}\max\limits_{1\leq k \leq n}\biggl|\sum_{j=1}^{k} r_{j}(\omega)\biggr|\, d\omega\leq c_2 \sqrt{n}.
 \]Hence, $\|t\|_{\mathrm{QC}} \leq c_3 \sqrt{n}$. Theorem \ref{S5:Th.Oskolkov} is proved.

 Oskolkov's example was generalized in \cite{Radomskii.Vestnik}.
 \begin{theorem}[Radomskii \cite{Radomskii.Vestnik}]\label{T.5.4.QC}
 Let $\{a_{n}\}_{n=1}^{\infty}$ be a sequence of real numbers such that
 \[
|a_{2^{j}-1}|+\sum_{k=2^{j-1}}^{2^{j}-2}|a_{k}-a_{k+1}| \leq \frac{1}{2^{j}},\quad j=1,2,\ldots
\]\textup{(}by definition, we put $\sum_{1}^{0}:=0$\textup{)}. Then
\[
\biggl\|\sum_{k=1}^{2^{n}-1}a_k \cos kx\biggr\|_{\mathrm{QC}}\leq c\sqrt{n},\quad n=1,2,\ldots,
\]where $c>0$ is an absolute constant.
  \end{theorem}

  \begin{corollary}[\cite{Radomskii.Vestnik}]\label{Cor.cos.QC}
  Let $N$ be an integer with $N\geq 2$. Then
  \begin{equation}\label{Cor.5.1.INEQ.cos}
  c_{1}\sqrt{\ln N}\leq \biggl\|\sum_{n=1}^{N}\frac{\cos nx}{n}\biggr\|_{\mathrm{QC}}\leq
  c_{2}\sqrt{\ln N},
  \end{equation}where $c_{1}$ and $c_{2}$ are positive absolute constants.
  \end{corollary}

  \textsc{Proof of Corollary \ref{Cor.cos.QC}.} We may assume without loss of generality that $N=2^{k}-1$, and put
  \[
  f(x)=\sum_{n=1}^{N}\frac{\cos nx}{n}.
  \] Applying Proposition \ref{S4:Rademacher.Prop} (3, ii), we have
  \begin{align*}
  \|f\|_{\mathrm{QC}}&=\int_{0}^{1}\biggl\|\sum_{j=1}^{k} r_{j}(\omega) \delta_{j}(f,\cdot)\biggr\|_{\infty}\,d\omega
  \geq \biggl\|\int_{0}^{1}\biggl|\sum_{j=1}^{k}r_{j}(\omega)\delta_{j}(f,x)\biggr|\,d\omega\biggr\|_{\infty}\geq\\
  &\geq c_{1}\biggl\| \biggl(\sum_{j=1}^{k}\delta_{j}(f,x)^{2}\biggr)^{1/2}\biggr\|_{\infty}\geq
  c_{1} \biggl(\sum_{j=1}^{k}\delta_{j}(f,0)^{2}\biggr)^{1/2},
  \end{align*}where $c_{1}$ is a positive absolute constant. Since
  \[
  \delta_{j}(f,0)=\sum_{n=2^{j-1}}^{2^{j}-1}\frac{1}{n}>\frac{1}{2^{j}}\,2^{j-1}=\frac{1}{2},\quad j=1,\ldots, k,
  \] the first inequality in \eqref{Cor.5.1.INEQ.cos} is proved. Applying Theorem \ref{T.5.4.QC} with $a_n=1/(2n),$ $n= 1, 2,\ldots,$ we obtain the second inequality in \eqref{Cor.5.1.INEQ.cos}. Corollary \ref{Cor.cos.QC} is proved.

  The following question was open for a long time: do there exist subspaces $L_{n}\subset E_{n},$ $n=1, 2,\ldots,$ such that $\dim L_{n}\geq \alpha \dim E_{n},$ $n=1, 2, \ldots,$ and
  \[
  A\|t\|_{\infty}\leq \|t\|_{\mathrm{QC}}\leq B\|t\|_{\infty}
  \]for any $n$ and $t\in L_{1}\oplus\cdots \oplus L_{n},$ where $\alpha\in (0,1),$ $A>0,$ $B>0$ are absolute constants? It was shown in \cite{Radomskii.MatSb} that the answer is negative. In \cite{Radomskii.Izv} we strengthened a result in \cite{Radomskii.MatSb}.
  \begin{theorem}[Radomskii \cite{Radomskii.Izv}]\label{S5:TH.RAD.QC}
   Let $\varepsilon$ be a real number with $\varepsilon\in (0,1)$, $\{k_n\}_{n=1}^{\infty}$ be a sequence of positive integers with  $k_1< \cdots< k_n<\cdots$, and $L_n$ be a subspace of $E_{k_n}$ such that $\dim L_n \geq \varepsilon \dim E_{k_n},$ $n=1, 2, \ldots .$ Then
   \[
    \sup\limits_{t\in L_1 \oplus\cdots \oplus L_n}\frac{\|t\|_{\mathrm{QC}}}{\|t\|_{\infty}}\geq \frac{d \varepsilon n}{a \varepsilon^{-1/2}n^{1/2} + b \varepsilon^{-3/2}}\geq \gamma \varepsilon^{5/2} \sqrt{n},\quad n=1, 2, \ldots,
   \]where $a$, $b$, $d$, and $\gamma$ are positive absolute constants.
  \end{theorem}

From Theorem \ref{S5:TH.RAD.QC} we obtain
\begin{corollary}[Radomskii \cite{Radomskii.Izv}]
Suppose that sequences $\{\tau_n\}_{n=1}^{\infty}$ and $\{\varepsilon_n\}_{n=1}^{\infty}$ of real numbers satisfy the conditions\\
\textup{i)} $\tau_n \geq 1$, $\tau_n \leq \tau_{n+1}$ \textup{(}$n=1, 2, \ldots$\textup{)} and $\tau_n \to +\infty$ as $n\to +\infty$\textup{;}\\
\textup{ii)} $\varepsilon_n \in (0,1)$, $\varepsilon_n \geq \varepsilon_{n+1}$ \textup{(}$n=1, 2, \ldots$\textup{)} and, starting at some $N$,
\[
\varepsilon_n = \frac{\tau_n}{n^{1/3}}.
\]Let $L_n$ be a subspace of $E_n$ such that
\[
\dim L_n \geq \varepsilon_n \cdot 2^n = \varepsilon_n \dim E_n,\quad n=1, 2, \ldots .
\]Then
\[
\sup\limits_{t\in L_1 \oplus\cdots \oplus L_n}\frac{\|t\|_{\mathrm{QC}}}{\|t\|_{\infty}}\geq \rho\cdot \tau_n^{3/2},\quad n\geq N,
\]where $\rho$ is a positive absolute constant. In particular,
\[
\lim\limits_{n \to \infty} \sup\limits_{t\in L_1 \oplus\cdots \oplus L_n}\frac{\|t\|_{\mathrm{QC}}}{\|t\|_{\infty}}= +\infty .
\]
\end{corollary}

\section{Sidon-type inequalities for discrete orthonormal systems}

In this section we consider real-valued functions $f$ with domain $(0,1)$ and use the following notation. Given $f\in L^{p}(0,1)$, we put
\begin{align*}
\|f\|_{p}&=\left(\int_{0}^{1}|f(x)|^{p}\, dx\right)^{1/p}\ \ \text{for $1\leq p <\infty$},\\
\|f\|_{\infty}&= \text{ess} \sup_{[0,1]} |f(x)|\ \ \ \,\ \ \ \ \ \ \ \,\text{for $p=\infty$}.
\end{align*}
For a system of functions $\Phi=\{\varphi_{k}(x)\}_{k=1}^{\infty},$ we let $\Phi(N),\ N\in\mathbb{N},$ denote the set of functions $p(x)$ of the form
\[
p(x)=\sum_{k=1}^{N} a_{k}\varphi_{k}(x),\quad a_{k}\in \mathbb{R},\ \  k=1,\ldots,N.
\]We recall that a system of functions $\{\varphi_{n}(x)\}_{n=1}^{\infty}\subset L^{2}(0,1)$ is \emph{orthonormal} if \[
\int_{0}^{1}\varphi_{n}(x)\varphi_{m}(x)\,dx=\begin{cases}
                                               1\ &\text{if $n=m$,}\\
                                               0\ &\text{if $n\neq m$,}
\end{cases}\qquad n, m = 1,2,\ldots.
\]The following result is the analog of a Sidon-type inequality for discrete orthonormal systems.
\begin{theorem}[Radomskii \cite{Radomskii.Discrete}]\label{S7:Th.DISCRETE}
 Let $\{m_{k}\}_{k=1}^{\infty}$ be a sequence of positive integers with $1=m_1<m_2<\cdots<m_k<\cdots$, and let $\Phi=\{\varphi_{k}(x)\}_{k=1}^{\infty}$ be an orthonormal system on $[0,1]$ with $\varphi_{1}(x)\equiv 1$ satisfying the following conditions:

 \textup{1)} for all $ k\geq 1,$ $ 1\leq j\leq m_k,$
\[
\varphi_{j}(x)=\textup{const}\quad \text{for }\ x\in\left(\frac{i-1}{m_k}, \frac{i}{m_k}\right),\ i=1,\ldots, m_k;
\]

\textup{2)} for all $ k\geq 2,$ $n> m_k,$
\[
\int_{0}^{1}\varphi_{n}^{2}(x)\varphi_{j}(x)\,dx=0\quad \text{for\ \  }2\leq j\leq m_k;
 \]

  \textup{3)} $\sup_{n\geq1} \|\varphi_{n}\|_{\infty}=M<+\infty$.

   Let $\{n_k\}_{k=1}^{\infty}$ be a sequence of positive integers with $n_{1}=1$ and $m_{k-1}<n_k\leq m_k,$ $k\geq2$. Let $N$ and $l$ be positive integers with $N\geq l+1$. Let $p_{k}\in \Phi(m_l),$ $k=l+1,\ldots,N,$ and let
   \[
f(x)=\sum_{k=l+1}^{N}p_{k}(x)\varphi_{n_{k}}(x).
\] Then
\[
\|f\|_{\infty}\geq\frac{1}{M}\sum_{k=l+1}^{N}\|p_{k}\|_{1}.
\]
\end{theorem}

\begin{corollary}[\cite{Radomskii.Discrete}]
If the assumptions of the theorem hold, then
\[
\frac{1}{M}\sum_{k=2}^{N}|a_k|\leq \biggl\|\sum_{k=2}^{N}a_{k}\varphi_{n_{k}}\biggr\|_{\infty}\leq M\sum_{k=2}^{N}|a_{k}|,\qquad N=2,3,\ldots.
\]
\end{corollary}
We recall the definition of the Walsh system $\{w_{n}(x)\}_{n=0}^{\infty}$. We write a positive integer $n$ in binary
\begin{gather*}
n=\sum_{k=0}^{\infty}\theta_{k}(n)2^k=\sum_{k=0}^{s(n)}\theta_{k}(n)2^k,\ \text{where $\theta_k=0$ or $1$},\\
\theta_{s(n)}(n)=1,\quad s(n)=[\log_{2}n].
\end{gather*}
\begin{definition}
The \emph{Walsh system} $W=\{w_{n}(x)\}_{n=0}^{\infty},$ $x\in [0,1],$ is a system of functions such that $w_{0}(x)\equiv 1$ and
\[
w_{n}(x)=\prod_{k=0}^{\infty}\bigl(r_{k+1}(x)\bigr)^{\theta_{k}(n)}= r_{s(n)+1}(x)\prod_{k=0}^{s(n)-1}\bigl(r_{k+1}(x)\bigr)^{\theta_{k}(n)},\ n=1, 2,\ldots,
\]where $r_k(x)$ is the $k$th Rademacher function.
\end{definition}

The Walsh system is a complete orthonormal system (see, for example, \cite{KashSaak}). We put $\varphi_{k}(x)=w_{k-1}(x)$, $k= 1, 2, \ldots .$ It is easy to verify that the system of functions $\Phi=\{\varphi_{k}(x)\}_{k=1}^{\infty}$ satisfies the assumptions of Theorem \ref{S7:Th.DISCRETE} with $m_{k}=2^{k-1}$ and $M=1$. For an integer $m\geq 0$, we let $W(m)$ denote the set of polynomials $p(x)$ in the Walsh system of the form
\[
p(x)=\sum_{k=0}^{m}a_k w_k (x),\quad a_k\in\mathbb{R},\ \ k=0,\ldots, m.
\]Applying Theorem \ref{S7:Th.DISCRETE}, we obtain the following result.
\begin{corollary}[Radomskii \cite{Radomskii.Discrete}]
Let $\{n_k\}_{k=1}^{\infty}$ be a sequence of positive integers such that $2^{k-1}\leq n_k<2^k,$ $k=1,2,\ldots.$ Let $m$ and $l$ be integers with $l\geq 0$ and $m\geq l+1$. Let $p_k\in W(2^l - 1),$ $k=l+1,\ldots,m,$ and let
\[
f(x)=\sum_{k=l+1}^{m} p_k(x)w_{n_k}(x).
\]Then
\[
\|f\|_{\infty}\geq \sum_{k=l+1}^{m}\|p_k\|_1.
\]
\end{corollary}

We recall that by an \emph{Hadamard matrix} of order $n$ we mean an $n\times n$ matrix $H=(h_{ij})_{i,j=1}^{n}$ with elements $h_{i j}\in \{-1, +1\}$ and pairwise orthogonal rows (therefore, the columns of this matrix are also pairwise orthogonal). Let $p$ be a prime with $p\geq 3,$ $r$ be a positive integer, and $\alpha = p^{r}+1$. Using the fact that there exist Hadamard matrices of order $\alpha\cdot 2^{k},\ k\geq 1$ (see, for example, \cite{Hall}), we can construct orthonormal systems $\{\varphi_{n}(x)\}_{n=1}^{\infty}$ satisfying the assumptions of Theorem \ref{S7:Th.DISCRETE} with $M=1,$ $m_1 = 1,$ and $m_k=\alpha\cdot 2^{k-1},\ k\geq 2$; in particular, these orthonormal systems are distinct from the Walsh system.

\end{document}